\documentclass[12pt,leqno]{article}

\usepackage{amsmath,amssymb,amsthm}

\makeatletter

\newtheorem{Theorem}{Theorem}

\newtheorem{Lemma}[Theorem]{Lemma}
\newtheorem{Corollary}[Theorem]{Corollary}

\theoremstyle{definition}

\theoremstyle{remark}

\def\({{\rm (}}
\def\){{\rm )}}

\let\Mathrm\operator@font
\let\Cal\mathcal

\def\standop#1{\mathop{\Mathrm #1}\nolimits}
\def\difstop#1#2{\expandafter\def\csname #1\endcsname{\standop{#2}}}
\def\defstop#1{\difstop{#1}{#1}}

\defstop{codim}
\defstop{id}
\defstop{Sym}
\defstop{Hom}
\defstop{End}
\defstop{rank}
\defstop{Spec}
\def\GL{{\sl{GL}}}

\defstop{Ker}
\defstop{Min}
\difstop{Image}{Im}
\difstop{tdeg}{trans.deg}

\def\tor{_{\Mathrm{tor}}}
\defstop{length}
\defstop{supp}
\defstop{Proj}
\defstop{Flat}
\defstop{Pfaff}
\defstop{Pf}
\defstop{Alt}
\defstop{Mat}
\defstop{Ext}
\defstop{Tor}
\defstop{Cobar}
\defstop{Coker}
\difstop{height}{ht}
\difstop{Image}{Im}

\def\specialarrow#1{\setbox\z@=\hbox{$\m@th
 \mathop{\vphantom{\rightarrow}}\limits^{\hspace{.5ex}{#1}\hspace
{.8ex}}$}\mathrel{\ifdim\wd\z@<1.2em\dimen\tw@
1.2em\else\dimen\tw@\wd\z@\fi\copy\z@\kern-\wd\z@\hbox to\dimen\tw@
{\rightarrowfill}}}

\def\sdarrow#1{\downarrow\hbox to 0pt{\scriptsize$#1$\hss}}
\def\suarrow#1{\uparrow\hbox to 0pt{\scriptsize$#1$\hss}}
\def\ssearrow#1{\searrow\hbox to 0pt{\scriptsize$#1$\hss}}

\def\ext{{\textstyle\bigwedge}}

\def\section{\@startsection{section}{1}{\z@ }%
{-3.5ex plus -1ex minus -.2ex}{2.3ex plus .2ex}{\bf }}

\long\def\refname{\par\kern -3ex
\begin{center}\rm R\sc{eferences}\end{center}\par\kern 
-2ex}

\def\@seccntformat#1{\csname the#1\endcsname.\quad}

\def\@@@sect#1#2#3#4#5#6[#7]#8{%
   \ifnum #2>\c@secnumdepth 
      \def \@svsec {}\else \refstepcounter {#1}%
      \def\@svsec{}
   \fi 
   \@tempskipa #5\relax 
   \ifdim \@tempskipa >\z@ 
     \begingroup #6\relax \@hangfrom {\hskip #3\relax 
     \@svsec}{\interlinepenalty \@M #8\par }\endgroup 
     \csname #1mark\endcsname {#7}
   \else 
   \def \@svsechd {#6\hskip #3\@svsec #8\csname #1mark\endcsname {#7}}
   \fi \@xsect {#5}}

\def\@@@startsection#1#2#3#4#5#6{%
 \if@noskipsec \leavevmode \fi \par \@tempskipa #4\relax \@afterindenttrue 
 \ifdim \@tempskipa <\z@ \@tempskipa -\@tempskipa \@afterindentfalse 
 \fi \if@nobreak \everypar {}\else \addpenalty {\@secpenalty }\addvspace 
  {\@tempskipa }\fi \@ifstar {\@ssect {#3}{#4}{#5}{#6}}{\@dblarg 
  {\@@@sect {#1}{#2}{#3}{#4}{#5}{#6}}}}

\def\theparagraph{\thesection.\arabic{paragraph}}
\def\aparagraph{\@@@startsection{paragraph}{2}{\z@ }%
              {1.75ex plus .2ex minus .15ex}{-1em}{\bf(\theparagraph) } }
\def\paragraph{\@@@startsection{paragraph}{2}{\z@ }%
              {1.75ex plus .2ex minus .15ex}{-1em}{}{\bf(\theparagraph)} }

\let\c@Theorem\c@paragraph

\title{Base change of invariant subrings}

\author{M{\sc itsuyasu} H{\sc ashimoto}}
\date{\normalsize
Graduate School of Mathematics, Nagoya University\\
Chikusa-ku,  Nagoya 464--8602 JAPAN\\
\small \tt hasimoto@math.nagoya-u.ac.jp}

\begin{document}

\maketitle

\begin{abstract}
Let $R$ be a Dedekind domain, $G$ an affine flat $R$-group scheme, and $B$
a flat $R$-algebra on which $G$ acts.
Let $A\rightarrow B^G$ be an $R$-algebra map.
Assume that $A$ is Noetherian.
We show that if the induced map $K\otimes A\rightarrow(K\otimes B)^{K\otimes G}
$ is an isomorphism for any algebraically closed field $K$ which is 
an $R$-algebra, then
$S\otimes A\rightarrow (S\otimes B)^{S\otimes G}$ is an isomorphism 
for any $R$-algebra $S$.
\end{abstract}

\section{Introduction}

In this paper, we prove the following.

\begin{Theorem}\label{main.thm}
Let $R$ be a Dedekind domain, $G$ an affine flat $R$-group scheme, 
and $M$ an $R$-flat $G$-module. 
Let $A$ be a Noetherian $R$-algebra, and $V$ a finitely generated 
$A$-module.
Let $\varphi\colon V\rightarrow M^G$ be an $R$-linear map.
If the induced map $\varphi_K\colon K\otimes V\rightarrow (K\otimes M)^{K
\otimes G}$ is an isomorphism for any algebraically closed field 
$K$ which is an $R$-algebra, then the canonical map
$\varphi_S\colon S\otimes V\rightarrow (S\otimes M)^{S\otimes G}$ is
an isomorphism for any $R$-algebra $S$.
\end{Theorem}

As a corollary, we have the following.

\begin{Corollary}
Let $R$ be a Dedekind domain, $G$ an affine flat $R$-group scheme,
and $B$ a flat $R$-algebra on which $G$-acts.
Let $A$ be a Noetherian $R$-algebra, and 
$\varphi\colon A\rightarrow B^G$ an $R$-algebra map.
If the induced map $\varphi_K\colon K\otimes A\rightarrow (K\otimes B)^{K
\otimes G}$ is an isomorphism for any algebraically closed field
$K$ which is an $R$-algebra, then the canonical map
$\varphi_S\colon S\otimes A\rightarrow(S\otimes B)^{S\otimes G}$ is
an isomorphism for any $R$-algebra $S$.
\end{Corollary}

So we may work only over algebraically closed field instead of general 
commutative ring, once we know that the action and the candidate of the
generator and the relation of the invariant subring are given over a
Dedekind domain (e.g., $\Bbb Z$), and the group scheme in problem is flat
over the Dedekind domain.

De Concini and Procesi \cite{DP} calculated the invariant subrings for
several important group scheme actions over an arbitrary commutative ring.
In \cite {Hashimoto2}, a simple proof (for the action of the general linear
group and the symplectic group) utilizing a geometric argument over a
field is given.
In order to reduce the case of general base ring to the case of
base field, the knowledge of good filtrations is utilized in 
\cite{Hashimoto2}, but this was a completely general theory as above,
since we know that the general linear group and the symplectic group 
are flat over $\Bbb Z$.

In section~2, we prove the theorem above.
In section~3, we give an example of applications.

\section{The proof of the main theorem}

Let $R$ be a commutative ring, and $G$ a flat $R$-group scheme.
Let $C$ be the coordinate ring $R[G]$ of $G$.
It is an $R$-flat commutative $R$-Hopf algebra.
A $G$-module is nothing but a right $C$-comodule, 
see \cite[Chapter~2]{Jantzen}.
For a $G$-module $M$, $M^G=\{m\in M\mid \omega(m)=m\otimes1\}$, where
$\omega\colon M\rightarrow M\otimes C$ is the coaction.
By means of the natural inclusion $\Hom_G(R,M)\hookrightarrow 
\Hom_R(R,M)=M$, the $R$-module $\Hom_G(R,M)$ is identified with
$M^G$, where $R$ is equipped with the trivial $G$-module structure.

In general, an $R$-module $V$ is considered as a trivial $G$-module.
So, for a $G$-module $M$ and an $R$-module $V$, $V\otimes M$ is a $G$-module
with the coaction
\[
1_V\otimes\omega_M\colon V\otimes M\rightarrow V\otimes M\otimes C
\]
where $\omega_M$ is the coaction of $M$.

The category of $G$-modules is abelian, with enough injectives, see
\cite[Lemma~I.3.3.3]{Hashimoto} and \cite[Lemma~I.3.5.9]{Hashimoto}.
For a $G$-module $M$ and an $R$-algebra $S$, the right 
$S\otimes C$-comodule structure of $S\otimes M$ is given by the composite
\[
S\otimes M\xrightarrow{1_S\otimes\omega}S\otimes M\otimes C
\xrightarrow{\alpha}(S\otimes M)\otimes_S(S\otimes C),
\]
where $\alpha$ is the isomorphism given by 
$\alpha(s\otimes m\otimes c)=(s\otimes m)\otimes(1\otimes c)$.
In particular, $(S\otimes M)^{S\otimes G}=(S\otimes M)^G$.
So $(S\otimes M)^G$ is an $S$-module.

Let $\varphi\colon V\rightarrow M^G$ be an $R$-linear map.
Then we define $\varphi_S\colon S\otimes V\rightarrow (S\otimes M)^G$ by
$\varphi_S(s\otimes v)=s\otimes\varphi(v)$.
For an $R$-algebra map $S\rightarrow S'$, we define
$\rho_{S',S}\colon S'\otimes_S (S\otimes M)^G\rightarrow (S'\otimes M)^G$ 
by
\[
\rho_{S'S}(s'\otimes(\sum_i s_i\otimes m_i))=\sum_i s's_i\otimes m_i.
\]
We denote $\rho_{S,R}\colon S\otimes M^G\rightarrow (S\otimes M)^G$ by 
$\rho_S$.
So $\rho_S(s\otimes m)=s\otimes m$ for $s\in S$ and $m\in M^G$.
Note that $\varphi_S$ is the composite
\[
S\otimes V\xrightarrow{1\otimes \varphi}S\otimes M^G
\xrightarrow{\rho_S}(S\otimes M)^G.
\]

For a $G$-module $M$, we denote $\Ext^i_G(R,M)$ by $H^i(G,M)$, and
call it the $i$th $G$-cohomology of $M$.
In particular, $H^0(G,M)=M^G$.

Let $M$ be a $G$-module.
Then by \cite[Lemma~I.3.6.16]{Hashimoto}, $H^i(G,M)\cong H^i(\Cobar_C(
M,R))$, where $\Bbb F(M):=\Cobar_C(M,R)$ is the complex
\[
M
\xrightarrow{\delta^0}
M\otimes C
\xrightarrow{\delta^1}
M\otimes C\otimes C
\xrightarrow{\delta^2}
\cdots
\]
whose boundary map is given by
\[
\delta^n=(-1)^{n+1}\omega_M\otimes 1_{C^{\otimes n}}
+\sum_{i=0}^{n-1}(-1)^{n-i}1_{M}\otimes 1_{C^{\otimes i}}\otimes \Delta_C
\otimes 1_{C^{\otimes n-i-1}}+1_M\otimes 1_{C^{\otimes n}}\otimes u,
\]
where $\Delta_C\colon C\rightarrow C\otimes C$ is the coproduct, and
$u\colon R\rightarrow C$ is the unit map.
By definition, for an $R$-module $V$, $\Bbb F(V\otimes M)\cong V\otimes
\Bbb F(M)$.
If $M$ is $R$-flat, then $\Bbb F(M)$ is an $R$-flat complex.
By the universal coefficient theorem \cite[Lemma~III.2.1.2]{Hashimoto} 
and its proof, we have the following.

\begin{Lemma}\label{universal.thm}
If $R$ is a Dedekind domain and $M$ is an $R$-flat $G$-module, there is 
an exact sequence
\[
0
\rightarrow
S\otimes M^G\xrightarrow{\rho_S}
(S\otimes M)^G\rightarrow
\Tor_1^R(S,H^1(G,M))\rightarrow0.
\]
\end{Lemma}

\proof[Proof of Theorem~\ref{main.thm}]
Let $R$, $G$, $A$, $V$, and $M$ be as in the theorem.

First, we prove the theorem for the case where $S$ is a field.
Let $K$ be the algebraic closure of $S$.
Taking the tensor product of $\varphi_S\colon 
S\otimes V\rightarrow(S\otimes M)^G$ with $K$ over $S$, 
we get $1\otimes\varphi_S\colon K\otimes V\rightarrow
K\otimes_S(S\otimes M)^G$.
As $K$ is faithfully flat over $S$, it suffices to show that 
this map is an isomorphism.
The composite
\[
K\otimes V\xrightarrow{1\otimes \varphi_S}
K\otimes_S(S\otimes M)^G
\xrightarrow{\rho_{K,S}}
(K\otimes M)^G
\]
is $\varphi_K$, which is an isomorphism.
Since $K$ is $S$-flat, $\rho_{K,S}$ is an isomorphism by
Lemma~\ref{universal.thm}.
So $1\otimes \varphi_S$ is an isomorphism as desired, and the
theorem is true for the case that $S$ is a field.

Next, we show that $H^1(G,M)$ is $R$-flat.
Since $R$ is Noetherian, it suffices to show that 
$\Tor^R_1(R/P,H^1(G,M))=0$ for any prime ideal $P$ of $R$.
Since $R$ is a one dimensional domain, it suffices to show that
$\Tor^R_1(R/\frak m,H^1(G,M))=0$ for any maximal ideal $\frak m$ of $R$.
On the other hand, $\varphi_{R/\frak m}$, which is the composite
\[
R/\frak m\otimes V
\xrightarrow{1\otimes \varphi}
R/\frak m\otimes M^G
\xrightarrow{\rho_{R/\frak m}}
(R/\frak m\otimes M)^G,
\]
is an isomorphism by the last paragraph.
So $\rho_{R/\frak m}$ is surjective.
By Lemma~\ref{universal.thm}, $\Tor^R_1(R/\frak m,H^1(G,M))=0$.
Hence $H^1(G,M)$ is $R$-flat, as desired.

Since $H^1(G,M)$ is $R$-flat, 
\[
\rho_S\colon S\otimes M^G\rightarrow(S\otimes M)^G
\]
is an isomorphism for any $R$-algebra $S$ by Lemma~\ref{universal.thm}.
Since the composite
\[
K\otimes V
\xrightarrow{1\otimes\varphi}
K\otimes M^G
\xrightarrow{\rho_K}
(K\otimes M)^G,
\]
which agrees with $\varphi_K$, is an isomorphism and $\rho_K$ is also
an isomorphism for any field $K$ which is an $R$-algebra, 
we have that $1\otimes\varphi\colon K\otimes V
\rightarrow K\otimes M^G$ is an isomorphism.

Next, we show that $V$ is $R$-flat.
First, we prove this for the case that $R$ is a DVR.
Let $t$ be a generator of the maximal ideal of $R$.
Since $V$ is a Noetherian $A$-module,
the torsion part $V\tor=\bigcup_{r\geq 0}(0:_{V}t^r)$ as an $R$-module
agrees with $(0:_{V}t^r)$ for some $r$.
Assume that $V\tor\neq 0$ for a contradiction.
Then $r\geq 1$.
We take $r$ as small as possible.
Take $a\in (0:t^r)\setminus(0:t^{r-1})$.
If $a\in tV$, then $a=ta'$ for some $a'\in V$.
Then $a'\in V\tor=(0:t^r)$.
So $t^{r-1}a=t^ra'=0$.
This contradicts the choice of $a$.
So $a\notin tV$.
Thus $1\otimes a\in R/tR\otimes_{R}V$ is nonzero.
Since $1\otimes \varphi\colon R/tR\otimes_{R}V\rightarrow R/tR
\otimes M^G$ is an isomorphism,
$1\otimes \varphi(a)\in R/tR\otimes M^G$ is nonzero.
This shows that $\varphi(a)\neq 0$ in $M^G$.
Since $M^G$ is a torsion free $R$-module, 
$\varphi(t^r a)=t^r\varphi(a)$ is nonzero.
This contradicts the assumption $t^ra=0$.
Hence $V$ is $R$-torsion free.
Since $R$ is a DVR, $V$ is $R$-flat.
Now consider the general case.
Let $\frak m$ be a maximal ideal of $R$.
Applying the discussion above to $R'=R_{\frak m}$, $A'=R'\otimes A$, 
$V'=R'\otimes V$ and $M'=R'\otimes M$, we have that $V_{\frak m}$ is
$R_{\frak m}$-flat for any $\frak m$.
This shows that $V$ is $R$-flat.

By \cite[Lemma~I.2.1.4]{Hashimoto}, $\varphi\colon V\rightarrow M^G$ is
injective, and $C:=\Coker\varphi$ is $R$-flat.
Since $K\otimes C=0$ for any field $K$ which is an $R$-algebra, we have
that $C=0$ by \cite[Corollary~I.2.1.6]{Hashimoto}.
Hence $\varphi$ is an isomorphism.

Let $S$ be any $R$-algebra.
The composite
\[
S\otimes V
\xrightarrow{1\otimes \varphi}
S\otimes M^G
\xrightarrow{\rho_S}
(S\otimes M)^G
\]
is an isomorphism, since $1\otimes \varphi$ and $\rho_S$ are.
This is what we wanted to prove.
\qed

\section{An application}

Let $R$ be a commutative ring.
For an $R$-scheme $Z$, we denote $\Gamma(Z,\Cal O_Z)$ by $R[Z]$.
For $v\geq 0$ and finite free $R$-modules $F$ and $G$, we denote by
$Y_v(F,G)$ the closed subscheme of $\Hom_R(F,G)$ consisting of $R$-linear
maps of rank at most $v$.
We denote the kernel of the map $R[\Hom_R(F,G)]\rightarrow R[Y_v(F,G)]$ by 
$I_{v+1}(F,G)$.
If $F$ and $G$ are of rank $f$ and $g$, respectively, then
$R[\Hom_R(F,G)]$ is identified with the polynomial ring $R[x_{ij}]_{1
\leq i\leq g,\;1\leq j\leq f}$ in $fg$ variables, and $I_{v+1}(F,G)$ is 
identified with the ideal of $R[x_{ij}]$ generated by the all 
$(v+1)$-minors of the matrix $(x_{ij})$.
Note that if $v\geq \min(f,g)$, then $Y_v(F,G)=\Hom_R(F,G)$, and
$I_{v+1}(F,G)=0$.

Let $m,n,r,s,t\in \Bbb Z_{\geq 0}$ such that
$s\leq m,r$ and $t\leq n,r$.
Set $u:=\min(s,t)$.
Let $V:=R^n$, $W:=R^m$, $E:=R^r$, $X:=Y_s(E,W)\times Y_t(V,E)$,
and $Y:=Y_u(V,W)$.
We define $\pi\colon X\rightarrow Y$ by $\pi(\varphi,\psi)=\varphi\circ
\psi$.
Let $G:=\GL(E)$ and $G':=\GL(W)\times\GL(V)$.
Then $G\times G'$ 
acts on $X$ by $(g,(g_1,g_2))\cdot(\varphi,\psi)=(g_1\varphi g^{-1},
g\psi g_2^{-1})$ for $g\in G$, $g_1\in\GL(W)$, $g_2\in \GL(V)$, 
$\varphi\in Y_s(E,W)$, and $\psi\in Y_t(V,E)$.
Letting $G\times G'$ act on $Y$ by 
$(g,(g_1,g_2))\cdot \rho=g_1\rho g_2^{-1}$ for $g\in G$, $g_1\in \GL(W)$, 
$g_2\in \GL(V)$, and $\rho\in Y$, 
the morphism $\pi$ is $G\times G'$-equivariant.
Note that $G$ acts on $Y$ trivially.

As an application of Theorem~\ref{main.thm}, we prove the following.

\begin{Theorem}
The morphism $\pi\colon X\rightarrow Y$ induces an isomorphism 
$\pi^\#\colon R[Y]\rightarrow R[X]^G$.
\end{Theorem}

By Theorem~\ref{main.thm}, we may assume that $R=K$ is an algebraically closed
field.

Let us recall some basic facts from representation theory.
A $G$-module $M$ is said to have good filtrations 
if $\Ext^1_G(\Delta_G(\lambda),M)
=0$ for any dominant weight $\lambda$, where $\Delta_G(\lambda)$ denotes
the Weyl module of highest weight $\lambda$, see \cite[(II.4.16)]{Jantzen}.

For a partition $\lambda=(\lambda_1,\ldots,\lambda_k)$ with $\lambda_1
\leq r$, the Schur module \cite{ABW} $L_\lambda E^*$ is a dual Weyl module.
In fact, $L_\lambda E^*\cong (\ext^r E^*)^{\otimes k}\otimes L_\mu E$, 
where $\mu=(r-\lambda_k,\ldots,r-\lambda_1)$.
By the Cauchy formula \cite{ABW}, $K[\Hom(E,W)]\cong \Sym(E\otimes W^*)$, 
$I_s(E,W)$, 
$K[Y_s(E,W)]$, $K[\Hom(V,E)]\cong \Sym(V\otimes E^*)$, $I_t(V,E)$, and
$K[Y_t(V,E)]$ have good filtrations as $G$-modules.
Since modules with good filtrations are closed under tensor products 
\cite{Wang}, \cite{Donkin}, 
\cite{Mathieu} and
extensions, the kernel $I$ of the canonical surjective map
$$
\rho\colon
K[\Hom(E,W)\times \Hom(V,E)]\rightarrow K[Y_s(E,W)\times Y_t(V,E)]
$$
has good filtrations, since there is a short exact sequence
\[
0\rightarrow
I_s(E,W)\otimes K[\Hom(V,E)]
\rightarrow I \rightarrow
K[Y_s(E,W)]\otimes I_t(V,E)
\rightarrow
0
.
\]
Hence $H^1(G,I)=0$.
It follows that $\rho$ induces a surjective map
\[
\rho^G\colon 
K[\Hom(E,W)\times \Hom(V,E)]^G\rightarrow 
K[Y_s(E,W)\times Y_t(V,E)]^G=K[X]^G.
\]
By the following theorem due to De Concini and Procesi \cite{DP},
$\pi^\#\colon K[Y]\rightarrow K[X]^G$ is surjective.

\begin{Theorem}
The composition
\[
\Hom(E,W)\times \Hom(V,E)\rightarrow Y_r(V,W)\qquad 
((\varphi,\psi)\mapsto \varphi\psi)
\]
induces
an isomorphism $K[Y_r(V,W)]\rightarrow K[\Hom(E,W)\times \Hom(V,E)]^G$.
\end{Theorem}

It remains to prove that $\pi^\#\colon K[Y]\rightarrow K[X]^G$ is injective.
As we know that $K[Y]$ is an integral domain (see e.g., \cite[(6.3)]{BV}), 
it suffices to show that $\pi$ is dominating.
By linear algebra,
for each $i$ such that $0\leq i\leq u$, the set of linear maps $V\rightarrow 
W$ of rank $i$ forms one $G'$-orbit.
Moreover, the $G'$-orbit of rank $u$ linear maps is dense in $Y$.
Since $\pi$ is $G'$-invariant, it suffices to show that $\pi(X)$ contains
at least one linear map of rank $u$.
But this is trivial.
\qed



\begin{thebibliography}{99}
\def\ji#1#2(#3)#4-#5.{\newblock{\em#1} {\bf#2} (#3), #4--#5.}
\def\GTM#1{Graduate Texts in Math. {\bf #1}, Springer}
\def\SLN#1{Lecture Notes in Math. {\bf #1}, Springer}

\bibitem{ABW} K. Akin, D. A. Buchsbaum and J. Weyman, Schur functors and Schur
complexes, 
{\em Adv. Math.} {\bf 44} (1982), 207--278.

\bibitem{BV} W. Bruns and U. Vetter, 
{\em Determinantal Rings}, \SLN{1327} (1988).

\bibitem{DP} C. De Concini and C. Procesi, A characteristic free approach
to invariant theory, {\em Adv. Math.} {\bf 21} (1976), 330--354.

\bibitem{Donkin} S. Donkin, {\em Rational Representations of 
Algebraic Groups}, \SLN{1140} (1985).

\bibitem{Hashimoto} M. Hashimoto, {\em Auslander-Buchweitz Approximations
of Equivariant Modules}, Cambridge (2000).

\bibitem{Hashimoto2} M. Hashimoto, Another proof of theorems of De Concini
and Procesi, preprint {\tt arXiv:math.AC/0408429 v1}.

\bibitem{Jantzen} J. C. Jantzen, {\em Representations of Algebraic Groups},
2nd edition, AMS (2003).

\bibitem{Mathieu} O. Mathieu, Filtrations of $G$-modules,
{\em Ann. Sci. \'Ecole Norm. Sup.} (4) {\bf 23} (1990), 625--644.

\bibitem{Wang} J. Wang, Sheaf cohomology on $G/B$ and tensor products of
Weyl modules, {\em J. Algebra} {\bf 77} (1982), 162--185.

\end{thebibliography}
\end{document}